\documentclass[reqno]{amsproc}
\usepackage{mathrsfs}
\usepackage{amsmath,amssymb,amsfonts}
\usepackage{mathtools}
\usepackage{bm}
\usepackage{upgreek}
\usepackage{hyperref}
\usepackage{cite}
\hypersetup{colorlinks=true,linkcolor=magenta,anchorcolor=green,citecolor=cyan,filecolor=black,menucolor=black,urlcolor=brown}

\theoremstyle{plain}
\newtheorem{theorem}{Theorem}[part]
\newtheorem{lemma}[theorem]{Lemma}
\newtheorem*{theorem*}{Theorem}

\theoremstyle{definition}
\newtheorem{definition}[theorem]{Definition}

\theoremstyle{remark}
\newtheorem{remark}[theorem]{Remark}

\numberwithin{section}{part}
\numberwithin{equation}{part}

\usepackage{soul}
\usepackage[dvipsnames]{xcolor}
\definecolor{bubblegum}{rgb}{0.99, 0.76, 0.8}

\usepackage[compat=1.1.0]{tikz-feynman}
\usepackage{paralist,tikz-cd}

\tikzset{
    partial ellipse/.style args={#1:#2:#3}{
        insert path={+ (#1:#3) arc (#1:#2:#3)}
    }
  }

\def\su{\circleddash}
  \usetikzlibrary{calc}

\newtheorem*{corollary*}{Corollary}

\def\IQ{\mathbb{Q}}
\def\IP{\mathbb{P}}
\def\IR{\mathbb{R}}

\def\cB{\mathscr{B}}

\def\cF{\mathscr{F}}

\def\cP{\mathscr{P}}
\def\cU{\mathscr{U}}
\def\cV{\mathscr{V}}
\def\cX{\mathscr{X}}

\newcommand{\Sol}{\mathcal{S}\!{\it o}\ell}

\def\Ree(#1){\Re\textrm{e}\left(#1\right)}
\def\Imm(#1){\Im\textrm{m}\left(#1\right)}

\def\Li#1(#2){\textrm{Li}_{#1}\left(#2\right)}

\begin{document}

\title[Twisted Picard-Fuchs Equations]{Picard-Fuchs Equations of Twisted Differential forms associated
  to Feynman Integrals}

\author{Pierre Vanhove}
\address{Institut de Physique Th\'eorique, Universit\'e Paris-Saclay, CEA, CNRS, F-91191 Gif-sur-Yvette Cedex, France}

\date{\today}

\thanks{This text is a contribution to the proceedings of the conference  Regulators V, 3-13 juin 2024 Department of Mathematics, University of Pisa, Italy.}
\begin{abstract}
Dimensionally or analytically regulated Feynman integrals lead to
relative twisted  period integrals.
We present a recent extension of the Griffiths-Dwork pole reduction
algorithm for deriving the D-module of differential operators acting
on the twisted differential forms from Feynman integrals.

We illustrate the application of this algorithm by providing twisted
Picard-Fuchs operators for hypergeometric, elliptic and
Calabi-Yau differential motives arising from families of Feynman integrals.

\end{abstract}

\maketitle
\newpage\tableofcontents\newpage

\part{Feynman integrals in Physics}

Feynman integrals are a cornerstone in understanding fundamental
interactions and the elementary building blocks of nature.  Scattering
amplitudes are used in particle physics to compare theoretical
predictions with experimental measurements in particle colliders
(see~\cite{Proceedings:2018jsb,Abreu:2022mfk,Craig:2022cef} for
instance), more  recently to gravitational wave
physics~\cite{Bjerrum-Bohr:2022blt,Kosower:2022yvp,Bjerrum-Bohr:2022ows},
or the evaluation of the
correlation functions of quantum fields at the end of inflation as they provide tools for analysing the formation of
structure in the Universe~\cite{Baumann:2022jpr,Benincasa:2022gtd}.
Their accurate calculation, whether analytically or numerically, is
needed for precision physics, but this  
remains a significant hurdle.

\medskip
  There are  growing evidences~\cite{Brown:2009ta,Bloch:2014qca,Bloch:2016izu,Bourjaily:2018ycu,Bourjaily:2019hmc,Bourjaily:2018yfy,Klemm:2019dbm,Bonisch:2020qmm,Bonisch:2021yfw,Bourjaily:2022bwx,Forum:2022lpz,Duhr:2022pch,Frellesvig:2023bbf,Pogel:2023zyd,Klemm:2024wtd,Driesse:2024feo,Frellesvig:2024rea, Duhr:2025lbz} that 
Feynman integrals needed for precision physics correspond to (relative) period integrals of (singular)
Calabi--Yau geometries.

Given a family of Feynman integrals attached to a Feynman graph,  it is desirable to answer the questions:
\begin{enumerate}
\item What is the class of functions to which belongs
a given  Feynman integral.
\item Determine the complete set of partial differential operators acting on
  a given family Feynman integrals.
\end{enumerate}

Feynman integrals satisfy several remarkable important properties:

\begin{enumerate}
  \item   They are D-finite functions, that is they satisfy a
    differentiable module of
  partial differential equations with respect to their
  parameters~\cite{Smirnov:2010hn,Lee:2013hzt}.
\item  Feynman integrals are relative period integrals of a variation
  of mixed Hodge structures~\cite{Bloch:2005bh,Brown:2009ta}.
  \end{enumerate}

These two properties set the question of analysing the nature of
Feynman integral in a clear Hodge theoretical framework. This approach
has been developed
in~\cite{Bloch:2005bh,Bloch:2008,Brown:2009ta,Bloch:2013tra,Weinzierl:2013fe,Bloch:2014qca,Brown:2015ylf,Brown:2015ztw,Bloch:2016izu,marcolli2019feynman,Vanhove:2014wqa,Bonisch:2021yfw,Doran:2023yzu}.

We recall, in section~\ref{sec:hodge}, how  Feynman integrals and  motives are
attached to a graph, and we make explicit the appearance of  twisted
differential in section. In section~\ref{part:differentialequations},  we explain the extension of the Griffths-Dwork reduction
for the specific twisted Feynman
differentials detailled in~\cite{delaCruz:2024xit}. In
section~\ref{sec:diffop}, we make explicit twisted differential operators
from Feynman integrals:
hypergeometric differential operator in section~\ref{sec:hyper},
hyperelliptic and elliptic differential operators in
section~\ref{sec:motives-two-loop}, and Calabi--Yau differential
operators in section~\ref{sec:CY}.

\part{Hodge structures of Feynman integrals}\label{sec:hodge}
\section{Feynman Graph polynomials}
\label{sec:graph-polynomials-1}

\begin{definition}
    A {\em Feynman} graph $\Gamma$ is a finite collection of vertices
    $V(\Gamma)$, edges $E(\Gamma)$, and half-edges $H(\Gamma)$
    satisfying the usual definitions; edges are adjacent to two
    vertices, and half-edges are adjacent to a single vertex, and
    allowing multiple edges between pairs of vertices.

    We let
    $e(\Gamma) = |E(\Gamma)|$ the number of edges. To each edge of
    $\Gamma$ we attach a mass variable $m_e \in \mathbb{R}$ and to
    each half-edge we attach  a momentum vector $p_h \in
    \mathbb{R}^{1,D-1}$ in the $D$-dimensional  Minkowski space equipped with a metric
    of  signature $(1,D-1)$.
To each
half edge of $\Gamma$ attach a vector $p_h \in \mathbb{C}^{D}$ subject to the so-called momentum conservation relation
\begin{equation}
\sum_{h \in H(\Gamma)} p_h = 0.
\end{equation}
 We 
assume that each
vertex of $\Gamma$ has a single outgoing half-edge. Therefore, one may view $\Gamma$ as a graph in
the usual sense, allowing multiple edges between vertices.
For
physical processes these vectors belong to of the
$D$-dimensional Minkowski space $\mathbb{R}^{1,D-1}$.
To simplify notation, we view momenta as being attached to vertices, and write $p_v$ instead of $p_h$. Furthermore, we consider only the completely massive case with $m_e^2>0$ and all external vectors are of non-zero norm
$p_v\cdot p_v\neq 0$. We take $m_e, p_v$ as having complex values.
The analytic properties of the Feynman integrals are studied by using
analytic continuation in the multi-dimensional complex plane spanned
by the number of independent scalar products $p_i\cdot p_j$ with
$1\leq i,j\leq |H(\Gamma)|$ and the masses $m_i^2$ with $1\leq i\leq e(\Gamma)$. 
\end{definition}

We associate to the graph $\Gamma$ two polynomials which  are defined as follows~\cite{nakanishi1971graph,Weinzierl:2022eaz}.
Let $\{ x_e \mid e \in e(\Gamma) \}$ be variables attached to all edges of
$\Gamma$. A spanning tree of $\Gamma$ is a subgraph ${\sf T}$ of $\Gamma$ which contains all vertices of $\Gamma$, and so that $b_1(\mathsf{T}) =0$ and $b_0(\mathsf{T})=1$. For each spanning tree $\mathsf{T}$ of $\Gamma$ we attach the
monomial $x^{\mathsf{T}} = \prod_{e\notin {\mathsf{T}}} x_e$. The {\em first Symanzik polynomial} is the polynomial
\begin{equation}
\cU = \sum_{\substack{ \text{Spanning} \\ \text{ trees of } \Gamma}} x^{\mathsf{T}}\,.
\end{equation}
A spanning $k$-forest of $\Gamma$ is a subgraph $\mathsf{F}$ of $\Gamma$ containing all vertices of $\Gamma$ and so that $h_1(\mathsf{F}) =0$ and $h_0(\mathsf{F}) = k$.  We attach the polynomial $x^{\mathsf{F}} =
\prod_{e \notin \mathsf{F}} x_i$ to each spanning 2-forest. A 2-forest
is a disjoint union of two sub-trees $\mathsf{F}=\mathsf{T}_1\cup \mathsf{T}_2$,  and we define $s_\mathsf{F} = \sum_{(v_1,v_2) \in \mathsf{F}=\mathsf{T}_1\cup \mathsf{T}_2} p_{v_1}\cdot p_{v_2}$ Where the ${}\cdot{}$-product is the
scalar product on $\mathbb {C}^{D}$. Then 
\begin{equation}\label{e:VFdef}
\cV(\vec{s}, \vec{m};D)  = \sum_{\substack{ \text{Spanning} \\ \text{ 2-forests
      of } \Gamma}}s_{\mathsf{F}} x^{\mathsf{F}},\quad \cF(\vec{s}, \vec{m};D) = \cU\times\left( \sum_{e \in E(\Gamma)} m_e^2 x_e\right) - \cV(\vec{s}, \vec{m};D) \,.
\end{equation}
The polynomial $\cF(\vec{s},\vec{m};D)$ is called the {\em second Symanzik
  polynomial} of $\Gamma$, depends 
the mass parameters and kinematic invariants, respectively:
\begin{equation}\label{e:kinematicparameters}
\vec{m}:=\left\{m_1^2,\dots,m_{e(\Gamma)}^2\right\}\in \mathbb
R_{>0}^{e(\Gamma)}, \qquad \vec s=\{p_i\cdot
p_j, i,j\in v(\Gamma)\}.
\end{equation}
When $|v(\Gamma)|>D$, not all the scalar
products are independent, and the number of independent variables
satisfies certain Gram determinant
conditions~\cite{Asribekov:1962tgp}.  The discriminant locus of $\cF$
depends on the linear relations between the
scalar products in $\vec s$.
This polynomial is a homogeneous polynomial of degree
$L+1$ in the variables $x_e$ for  $e \in e(\Gamma)$, where $L =
b_1(\Gamma)$. This $L$  is often called the {\em loop order} of
$\Gamma$. Henceforward, we will write instead $\cF$ to
simplify our notation.

\section{Feynman integrals in parametric representation}
\label{sec:feynm-param}

Using the   two polynomials associated to a graph $\Gamma$,  we define a family of  Feynman integrals~\cite{Itzykson:1980rh,Weinzierl:2022eaz}

\begin{equation}\label{e:GraphIZ1}
  I_\Gamma(\underline z;D,\underline \nu)=\int_{[0,+\infty[^{e(\Gamma)}}   \,
  {\cU^{\nu- (l+1){D\over2}}\over \cF^{\nu-l
      {D\over2}}}\,\delta\left(\sum_{i=1}^{e(\Gamma)} x_i-1\right)\prod_{i=1}^{e(\Gamma)} x_i^{\nu_i-1}  dx_i\,.  
\end{equation}
where we have set $\underline \nu:=(\nu_1,\dots,\nu_{e(\Gamma)})$ and collected
the kinematic factors (the internal masses $m_i$ and the independent
scalar products between the external momenta) into $\underline
z=(\vec s,\vec m)$.

Since  the coordinate scaling $(x_1,\dots,x_{e(\Gamma)})\to \lambda
(x_1,\cdots,x_{e(\Gamma)})$ leaves invariant the integrand and the domain of integration, we can rewrite this
integral as

\begin{equation}\label{e:GraphIZ2}
  I_\Gamma(\underline z;D,\underline \nu)=\int_{\Delta_{e(\Gamma)}}   \Omega_\Gamma^{D,\underline \nu}
\end{equation}
with
\begin{equation}
  \label{e:OmegaD}
  \Omega_\Gamma^{D,\underline \nu}:=  
  {\cU^{\nu- (l+1){D\over2}}\over \cF^{\nu-l {D\over2}}}\, \prod_{e\in
  e(\Gamma)} x_e^{\nu_e-1}\, \Omega_0
\end{equation}
with the differential $e(\Gamma)-1$-form

\begin{equation}\label{e:diffInt}
 \Omega_0 :=\sum_{j=1}^{e(\Gamma)} (-1)^{j-1} x_j\,
 dx_1\wedge\cdots\wedge \widehat{dx_j}\wedge\cdots \wedge dx_{e(\Gamma)},
\end{equation}
where $\widehat{dx_j}$ means that $dx_j$ is omitted.
The domain of integration is defined as
\begin{equation}\label{e:DefDomain}
  \Delta_{e(\Gamma)}:=\left\{[x_1,\dots,x_{e(\Gamma)}]\in \mathbb P^{e(\Gamma)-1}| x_i\in\mathbb R,  x_i\geq0\right\}.  
\end{equation}

\section{Mixed Hodge structures for Feynman graph integrals}
\label{sec:motive}

We define the vanishing loci  for the Symanzik
polynomials attached to the graph $\Gamma$:
\begin{equation}\label{e:VFVU}
    X_{\Gamma;D} =  \{  \cF(\vec{s},\vec{m};D)=0 | x_i\in\mathbb  P^{e(\Gamma)-1}(\IR)\};  \qquad Y_\Gamma =  \{  \cU=0 | x_i\in\mathbb
P^{e(\Gamma)-1}(\IR)\}\,.
\end{equation}
Notice that the vanishing
locus for $\cF$ depends on the space-time dimension $D$
through the linear relations between the external momenta.

The integrand of the Feynman integral~\eqref{e:GraphIZ2} is  a differential form representing a class of
$\mathrm{H}^{e(\Gamma)-1}(\mathbb{P}^{e(\Gamma)-1}-  Z_{\Gamma;D})$
where $Z_{\Gamma;D}$ is the singular locus of the integrand. We see that if $e(\Gamma) - {(L+1)D\over2}<0$, $Z_{\Gamma;D} = Y_{\Gamma}$ and that if and $e(\Gamma) - {LD\over2}>0$ then $Z_{\Gamma;D} = X_{\Gamma;D}$. If neither of these inequalities is satisfied, then $Z_{\Gamma;D} = X_{\Gamma;D} \cup Y_\Gamma$.

\medskip
Although the integrand $\Omega_\Gamma^{D,\underline \nu}$ is a closed form  such that $\eta\in
\mathrm{H}^{e(\Gamma)-1}(\mathbb{P}^{e(\Gamma)-1}-  Z_{\Gamma;D})$,  in general the domain
$\Delta_{e(\Gamma)}$ has a boundary and therefore its  homology class is not in $\mathrm{H}_{e(\Gamma)-1}(\mathbb{P}^{e(\Gamma)-1}-  Z_{\Gamma;D})$.  This difficulty is resolved by considering the relative cohomology~\cite{Bloch:2005bh,Brown:2015ylf}.

We need to consider a blow-up in $\IP^{e(\Gamma)-1}$ of linear space
$f:\cP\to \IP^{e(\Gamma)-1}$, such that all the vertices of $\Delta_{e(\Gamma)}$ lie in
$\cP\backslash\cX$ where $\cX$ is the strict transform of $ Z_{\Gamma;D}$.
Let  $\cB$ be the total inverse image of the coordinate simplex
$\{x_1x_2\cdots x_{e(\Gamma)}=0| [x_1,\dots,x_{e(\Gamma)}]\in\IP^{e(\Gamma)-1}\}$.

As been explained by Bloch, Esnault and Kreimer in~\cite{Bloch:2005bh}
all of this lead to the mixed Hodge structure associated to the Feynman graph
\begin{equation}
  \label{e:defMotive}
  M(\Gamma):=H^{e(\Gamma)-1}(\cP\backslash \cX, \cB\backslash
  \cB\cap\cX; \IQ)\,.
\end{equation}
%

\section{Regulated Feynman integrals and twisted differential forms}
\label{sec:diverg-epsil-expans}

The Feynman integral defined in~\eqref{e:GraphIZ2} is a function of
the parameters $D$ and $\underline \nu$. For integer values of $D$ and
the powers $\nu_i$ the integral can be divergent. There are the
ultaviolet or infrared divergences which have a special meaning in
quantum field theory. We refer to e.g.~\cite{Peskin:1995ev} for a physics based discussion.

One can prove that the Feynman integral is a meromorphic
function of $(D,\underline \nu)$  in $\mathbb C^{e(\Gamma)+1}$, with simple poles located on
affine hyperplanes defined by linear equations $c_0 D+\sum_{r=1}^{e(\Gamma)}
c_r \nu_r$ with integer coefficients
$(c_0,c_1,\dots,c_{e(\Gamma)})\in\mathbb Z^{e(\Gamma)+1}$. One can show that there is an open subset of $(D, \nu_1,
\dots , \nu_{e(\Gamma)}) \in \mathbb C^{e(\Gamma)+1}$ where the integral converges. The
(unique) value of the Feynman integral is defined by analytic
continuation. We refer to~\cite{Speer} for a throughouly
discussion.

\medskip

Using these properties there are two commonly used regularisation in
physics:
\begin{enumerate}[(1)]
  \item The dimensional regularisation where the analytic continuation is done
in the spacetime $D$ with fixed values for the indices $\underline \nu$ taken to be
integers.
\item The analytic continuation with a fixed integer value for the
spacetime dimension $D$, but the analytic continuation is done with
respect to (a subset) of the indices $\underline \nu$.
\end{enumerate}
In the following we combine both of these regulators, and we introduce the following notations
\begin{equation}
  I_\Gamma^{\epsilon,\kappa}(\underline z):=I_\Gamma(\underline z;
  2\delta-2\epsilon, \upnu_1+\mu_1 \kappa,\dots, \upnu_{e(\Gamma)}+\mu_{e(\Gamma)} \kappa)
\end{equation}
with $\delta$ an integer spacetime dimension of interest (say for instance
$\delta=2$), and $(\upnu_1,\dots,\upnu_{e(\Gamma)},\mu_1,\dots,\mu_{e(\Gamma)})\in\mathbb Z^{2 {e(\Gamma)}}$, and
$\epsilon$ and $\kappa$  are positive real numbers.
We set
\begin{equation}
    I_\Gamma^{\epsilon,\kappa}(\underline z)=\int_{\Delta_{e(\Gamma)}}
    \Omega_\Gamma^{\epsilon,\kappa},
\end{equation}
with
\begin{equation}\label{e:OmegaAnRegIntro}
  \Omega_\Gamma^{\epsilon,\kappa}= \omega_\Gamma^{\rm Rat}
  \times\left(\cU^{L+1}\over \cF^L\right)^{\epsilon} \, \prod_{i=1}^{e(\Gamma)}
  \left(x_i\cU\over \cF\right)^{\mu_i\kappa}\, dx_1\cdots dx_{e(\Gamma)},
  \end{equation}
with the rational function
\begin{equation}\label{e:OmegaRat}
  \omega_\Gamma^{\rm Rat}=   { \cU^{\upnu_1+\cdots+\upnu_{e(\Gamma)}-(L+1)\delta}  \over \cF ^{\upnu_1+\cdots+\upnu_{e(\Gamma)}-L\delta} }\, \prod_{i=1}^{e(\Gamma)} x_i^{\upnu_i-1} 
\end{equation}
We make a few remarks:

\begin{enumerate}[(1)]
\item When $\epsilon=\kappa=0$ we have that $\Omega_\Gamma^{0,0}=
  \omega_\Gamma^{\rm Rat}$ is a rational differential form.
  \item We remark that ${\cU^{L+1}/ \cF^L}$ and ${x_i\cU/
      \cF}$  are degree zero
    rational functions in $\mathbb P^{{e(\Gamma)}-1}$, therefore  when $\epsilon\neq0$ or
    $\kappa\neq0$ we have well-defined twisted
    differential forms.
    \item The twists in~\eqref{e:OmegaAnRegIntro} does not introduce
      new poles. Therefore,  the twisted differential has the
      \emph{same} singular locus as the rational differential form~\eqref{e:OmegaRat}. 
  \end{enumerate}

  Feynman integrals naturally lead to twisted differential
forms  of the kind
studied in~\cite{Aomoto1,Aomoto,Aomoto_1982,AomotoBook}, which has
been considered  in recent applications, e.g.~\cite{Mizera:2017rqa,Frellesvig:2019uqt,Mizera:2019ose}, for reducing the families of
Feynman integrals attached to a given graph onto a basis of integrals.

In the present
text will discuss a different approach aimed to derive the
differential operators acting on the Feynman
integrals~\cite{delaCruz:2024xit}.   Our approach uses that the twist
for Feynman integral is a rational function build from the graph polynomials.
\part{D-modules for twisted differential forms}\label{part:differentialequations}

An important property of Feynman integral is that they are
holonomic functions,  which means that they
satisfy finite order differential equation when differentiating with
respect to their physical parameters $\underline z:=\{\vec s,\vec m\}$.

Let us consider $r$ parameters from the set of internal masses and independent kinematics, $\underline z:=\{z_1,\dots,z_r\} \in \vec m \cup \vec s $. 
We seek differential operators annihilating the twisted differential form
$\Omega_\Gamma^{\epsilon,\kappa}$ in~\eqref{e:OmegaAnRegIntro} in cohomology
\begin{equation}\label{e:PFOmegaGeneric}
	\left(  \sum_{a_1=0}^{o_1} \cdots\sum_{ a_r=0}^{o_r} c_{a_1,\dots,a_r}(\vec m,\vec s,\epsilon,\kappa)\left(\partial\over \partial z_1\right)^{a_1}\cdots\left(\partial
	\over\partial z_r\right)^{a_r}  \right)\Omega_\Gamma^{\epsilon,\kappa}=d\beta^{\epsilon,\kappa}_\Gamma,
\end{equation}
where $c_{a_1,\dots,a_r}(\vec m,\vec s,\epsilon,\kappa)$ are rational functions of the physical
parameters, but they are independent of the edge variables
$x_1,\dots,x_{e(\Gamma)}$, and $o_1,\dots, o_r$ are some positive
integers. The inhomogeneous term $d\beta^{\epsilon,\kappa}_\Gamma$ is a total derivative in $x_i$'s where the only allowed poles are those already present in $\Omega_\Gamma^{\epsilon,\kappa}$~\cite{Lairez:2022zkj}.
Because the domain of integration~\eqref{e:DefDomain} of the Feynman integral does not
depend on the physical parameters, we then deduce
\begin{equation}\label{e:PDEintFeynman}
	\left( \sum_{a_1=0}^{o_1} \cdots\sum_{ a_r=0}^{o_r} c_{a_1,\dots,a_r}(\vec m,\vec s,\epsilon,\kappa)\left(\partial\over \partial z_1\right)^{a_1}\cdots\left(\partial
	\over\partial z_r\right)^{a_r}  \right)I^{\epsilon,\kappa}_\Gamma=\mathscr{S}^{\epsilon,\kappa}_\Gamma \, ,
\end{equation}
where $\mathscr{S}^{\epsilon,\kappa}_\Gamma$ is an inhomogeneous term obtained by
integrating $d\beta^{\epsilon,\kappa}_\Gamma$ over the boundary of the positive
orthant~\eqref{e:DefDomain}. This is a non-trivial task because one needs
to blow-up the intersections between the graph hypersurface and the
domain of integration, so the integral is well-defined~\cite{Bloch:2005bh,Brown:2009ta,Bloch:2016izu,muller2014picard}. For
instance, 
section~3.2 of~\cite{Bloch:2016izu} gives  a detailed derivation of the inhomogeneous term for
the two-loop sunset integral along these lines.

\medskip
If the integration is done over a cycle $\mathcal C$, like the one
defined  by the torus $\mathcal C_{\rm
  max}:=\{|x_1|=\cdots=|x_{e(\Gamma)}|=1\}$,  the resulting integral is
annihilated by the  action of the differential operator~\cite{Vanhove:2018mto}
\begin{equation}
	\left( \sum_{a_1=0}^{o_1}\cdots \sum_{ a_r=0}^{o_r} c_{a_1,\dots,a_r}(\vec m,\vec s,\epsilon,\kappa)\left(\partial\over \partial z_1\right)^{a_1}\cdots\left(\partial
	\over\partial z_r\right)^{a_r}  \right)\int_{\mathcal C_{\rm max} }\Omega^{\epsilon,\kappa}_\Gamma=0\,.
\end{equation}

The ideal  generated by these
differential operators is a differential module (or D-module).  Thus, the
differential equations we are seeking can be obtained by deriving  annihilators  of
$\Omega_\Gamma^{\epsilon,\kappa}$, i.e., partial differential operators that
annihilate  the integrand by acting on the physical parameter and the
edge variables.
An example of a system of partial differential equations for Feynman integrals is the Gel'fand-Kapranov-Zelevinski\u\i {} (GKZ) system, which provides
a D-module of differential operators acting on the toric generalisation
of the Feynman
integral~\cite{Vanhove:2018mto,delaCruz:2019skx,Klausen:2019hrg,Feng:2019bdx,Klemm:2019dbm,Ananthanarayan:2022ntm,Agostini:2022cgv,Matsubara-Heo:2023ylc,Munch:2022ouq}.
Because the graph polynomials in the expression for
$\Omega_\Gamma^{\epsilon,\kappa}$ in~\eqref{e:OmegaAnRegIntro} are not
generic polynomials the differential module acting on a
given Feynman integral is obtained after restricting the GKZ  D-module which
is a higher non-trivial
task~\cite{delaCruz:2019skx,Klausen:2021yrt,Chestnov:2023kww,Dlapa:2023cvx}
and its systematic implementation is 
still an open problem. In the following sections we present an
algorithmic procedure to derive the differential equations based on
an extension of the Griffiths-Dwork reduction for twisted differential forms.

A motivation is to have an  algorithm that applies to a large class of analytic
regularised Feynman integrals which  is missing 
for the commonly used programs in theoretical physics.
 This way we
can analyse how the twists parameters $\epsilon$, from space-time dimension,  and
$\kappa$, from the analytic regulator,  deform the minimal order of the
differential operators.

\section{Variation of mixed Hodge structure and ODEs}\label{s:vmhsode} 

If $\mathcal{H}_\mathbb{Q}$ is the local system underlying a variation
of mixed Hodge structure over a 1-dimensional base $M$, and ${\bm s}$
is a meromorphic section of $\mathcal{H}_\mathbb{Q}\otimes
\mathcal{O}_M$ then there is a minimal differential equation
$\mathscr{L}_{\bm s}$ annihilating the period functions attached to
${\bm s}$. 

\medskip

Selecting a parameter $t$ amongst the physical parameters $\vec m\cup
\vec s$ we consider a pencil of graph $\cF_\Gamma(t)$ and the
differential form
\begin{equation}
\Omega_{\Gamma}^{D,\underline\nu}(t) =\dfrac{\cU^{\sum_{i=1}^{e(\Gamma)}\nu_i -
    (L+1)D/2}}{(\cF_\Gamma(t))^{\sum_{i=1}^{e(\Gamma)}\nu_i - LD/2}} \, \prod_{i=1}^{
  e(\Gamma)} x_i^{\nu_i-1}\, \Omega_0,
\end{equation}
determines a section of $\mathcal{H}_{\Gamma;D}\otimes \mathcal{O}_M$
where $\mathcal{H}_{\Gamma;D}$ is a variation of mixed Hodge structure over an open subset $M$ of $\mathbb{A}^1_t$.

\begin{definition}
    Let ${\mathscr{L}}_{\Gamma}^{D}$ denote the minimal
    differential operator in $\mathbb{C}[M]\langle \partial_t\rangle$ which
    annihilates the form $\Omega_{\Gamma}^{D,\underline \nu}(t)$ in
    $\mathcal{H}_{\Gamma;D} \otimes \mathcal{O}_M$. 
\end{definition}

We recall  how an ordinary differential equation is associated  with a variation of mixed Hodge structure along with a holomorphic section of the underlying local system.  

\begin{definition}
A ($\mathbb{Q}$-)variation of mixed Hodge structure of weight $n$ consists of several pieces of data
\begin{enumerate}[(1)]
\item A $\mathbb{Q}$-local system $\mathcal{H}_\mathbb{Q}$ over a complex manifold $M$,
\item An increasing weight filtration by $\mathbb{Q}$-local systems $\mathcal{W}_0\subseteq \mathcal{W}_1\subseteq \dots \subseteq \mathcal{W}_{2n} = \mathcal{H}_\mathbb{Q}$,
\item A decreasing Hodge filtration $\mathcal{F}^n \subseteq \mathcal{F}^{n-1}\subseteq \dots \subseteq \mathcal{F}^0 = \mathcal{H}_\mathbb{C} = \mathcal{H}_\mathbb{Q}\otimes_{\underline{\mathbb{Q}}_M} \underline{\mathbb{C}}_M$,
\item A flat connection $\nabla : \mathcal{H}_\mathbb{C} \otimes \mathcal{O}_M \rightarrow \mathcal{H}_\mathbb{C} \otimes \Omega_M^1$ so that $\nabla(\mathcal{F}^i) \subseteq \mathcal{F}^{i-1}$,
\end{enumerate}
so that on each fibre $\mathcal{H}_{\mathbb{Q}},$ the data $(\mathcal{H}_{\mathbb{Q}},\mathcal{F}^\bullet_t,\mathcal{W}_{\bullet})$ is a mixed Hodge structure. 
\end{definition}

\noindent Given a local section ${\bm s}$ of $\mathcal{H}_\mathbb{C}\otimes \mathcal{O}_M$, and a local parameter $t$ on $M$, we can construct local (or multivalued) period functions
\begin{equation}\label{eq:pair}
{\bm \pi}_{\bm s}(t) = \langle {\bm s}, \gamma_t \rangle
\end{equation}
for a flat section $\gamma_t$ of $\mathcal{H}^\vee_{\mathbb{Q}}$. We will often take ${\bm s} = \Omega_{\Gamma}^{D,\underline \nu}(t)$ and let $\mathcal{H}^\vee_\mathbb{Q}$ is the homology bundle underlying the family of varieties $\mathbb{P}^{e(\Gamma)-1}-X_{\Gamma;D}(t)$, in which case the pairing is integration.

Given a variation of mixed Hodge structure, $(\mathcal{H}_\mathbb{Q},
\mathcal{W}_\bullet, \mathcal{F}^\bullet)$ over $M\subseteq
\mathbb{A}^1$ with Gauss--Manin connection $\nabla$, we have
differential operators
\begin{equation}
  \nabla_{\partial_t} :
  \mathcal{H}\otimes \mathcal{O}_M \rightarrow \mathcal{H}\otimes
  \mathcal{O}_M, [\omega]\mapsto
  \nabla([\omega])(\partial_t)
\end{equation}
where $\partial_t$ denotes the vector field corresponding to a choice of variable $t$. The pairing satisfies
\begin{equation}
\dfrac{d}{dt }\langle {\bm s}, \gamma_t\rangle = \langle \nabla_{\partial_t}(\omega), \gamma_t \rangle.
\end{equation}
Consequently, there is a minimal collection of elements $\{f_0(t),\dots, f_n(t)\}$ in the  $\mathbb{C}(t)$-vector space $\Gamma(\mathcal{H}\otimes \mathcal{O}_M)\otimes \mathbb{C}(t)$ so that 
\begin{equation}
\left[f_n(t)\nabla_{\partial_t}^n + f_{n-1}(t)\nabla_{\partial_t}^{n-1}  + \dots + f_1(t) \nabla_{\partial_t} + f_0(t)\right] {\bm s} = 0
\end{equation}
and thus there is a linear differential operator
\begin{equation}
\mathscr{L}_{\bm s} = f_n(t)\dfrac{d^n}{dt^n} + f_{n-1}(t)\dfrac{d^{(n-1)}}{dt^{(n-1)}}  + \dots + f_1(t) \dfrac{d}{dt} + f_0(t)
\end{equation}
whose solutions are the period functions ${\bm \pi}_{\bm s}(t)$. The local system $\mathcal{H}_\mathbb{Q}^\vee$ is equipped with a weight filtration $\mathcal{W}_{\bullet}^*$ dual to the weight filtration on $\mathcal{H}_{\Gamma;D}(t)$ determined by $\mathcal{W}_{i}^* = (\mathcal{W}_{-i-1})^\vee$. The pairing \eqref{eq:pair} induces a map from $\mathcal{H}_\mathbb{Q}^\vee$ to $\mathcal{O}_M$ whose image is $\Sol(\mathscr{L}_{\bm s})$, the local system of solutions of $\mathscr{L}_{\bm s}$. Therefore, $\mathcal{W}_i^*$ induces a filtration on $\Sol(\mathscr{L}_{\bm s})$.
\begin{lemma}[section~3.1 of~\cite{Doran:2023yzu}]\label{l:qotapp}
The local system $\Sol({\mathscr{L}}_{{\bm s}})$ is a quotient of the dual local system $\mathcal{H}_{\mathbb{Q}}^\vee$ by a sub-local system $\mathbb{K}_{\bm s}$. If $\bm{s} \in \mathcal{W}_i\otimes \mathcal{O}_M$ then $\mathcal{W}_i^* \subseteq \mathbb{K}_{\bm s}$.
\end{lemma}

\noindent We summarize the results of~\cite{Doran:2023yzu}
\begin{enumerate}[(1)]
    \item The local systems $\Sol(\mathscr{L}_{\Gamma;D})$ are quotients of $\mathcal{H}_{\Gamma;D}^\vee$.
    \item The filtration induced by $\mathcal{W}^*_\bullet$ corresponds to a factorisation of $\mathscr{L}_{\Gamma;D}$, however there may be factorisations of $\mathscr{L}_{\Gamma;D}$ which do not correspond to $\mathcal{W}^*_\bullet$.
    \item The monodromy representation of $\Sol(\mathscr{L}_{\Gamma;D})$ is upper triangular with diagonal blocks equal to the monodromy representations of the factors of $\mathscr{L}_{\Gamma;D}$.
\end{enumerate}

\section{Griffiths-Dwork reduction for twisted differential forms}\label{sec:griff-dwork-reduct}
We present the Griffiths-Dwork reduction for twisted differential
forms applied to the case of the differential form
$\Omega_\Gamma^{\epsilon,\kappa}$ defined
in~\eqref{e:OmegaAnRegIntro}.

Choosing $r$ variables amongst the kinematic parameters $\underline z:=\{z_1,\dots,z_r\} \in \vec m \cup \vec s $,
the    differentiation of $\Omega_\Gamma^{\epsilon,\kappa}$ 
leads  to 
\begin{equation}\label{e:PDE}
\sum_{\mathsf{a}=a_1+\cdots+a_r\atop
  a_i\geq0}\!\!\!\!\!\!\!\!c_{\underline a}(\vec m,\vec s;\epsilon,\kappa)\left(\partial\over \partial z_1\right)^{a_1}\cdots\left(\partial
  \over\partial z_r\right)^{a_r} \Omega_\Gamma^{\epsilon,\kappa}= \sum_{\mathsf{a}=a_1+\cdots+a_r\atop a_i\geq0}\!\!\!\!\!\!\!\!
{ 
  c_{\underline a}(\vec m,\vec s;\epsilon,\kappa)  P^{\underline a}(\underline x)\over \cF^{\mathsf {a}}}\, \Omega_\Gamma^{\epsilon,\kappa},
\end{equation}
where $\underline a=(a_1,\dots,a_r)\in\mathbb N^r$ and 
$  P^{\underline a}(\underline x)$ is a
  homogeneous polynomial of degree $(L+1)(a_1+\cdots+a_r)$ in
  the edge variables $\underline x$.
  The sum is over the differential operators of order
  $a_1\geq0 ,\dots, a_r\geq0$ and fixed total order $\mathsf{a}:=a_1+\cdots +a_r$. 
The
  pole order in the second Symanzik polynomial $\cF$ has increased by
  $\mathsf a$. 

  \medskip
We present an extension of the Griffths pole reduction~\cite{Griffith1,Griffith2} adapted  to include the twist
   factor in $\Omega_\Gamma^{\epsilon,\kappa}$:

   \begin{itemize}
   \item[\textbf{Step~1:}]   Reduction of polynomial
$P^{\underline a}(\underline x)$ in the Jacobian ideal of $\cF$,
$\text{Jac}( \cF_\Gamma):=\langle \vec\nabla \cF(\underline x)\rangle$
\begin{equation}\label{e:RedF}
	P^{\underline a}(\underline x) = \vec C^{\underline a}(\underline x)\cdot
	\vec\nabla   \cF \, ,
      \end{equation}
      with $	\vec\nabla   \cF :=\left(\partial_{x_1}
      \cF(\underline x),\dots, \partial_{x_{e(\Gamma)}}
      \cF(\underline x)\right)$.
    \item[\textbf{Step~2:}] Reduction of $ \vec C^{\underline a}(\underline x)$ in the
      Jacobian ideal of $\cU$, $\text{Jac}( \cF_\Gamma):=\langle
      \vec\nabla \cU(\underline x)\rangle$
      \begin{equation}\label{e:RedU}
\vec C^{\underline a}(\underline x)\cdot
	\vec\nabla   \cU = c^{\underline a}(\underline x)\, \cU\, ,
      \end{equation}
      with $	\vec\nabla   \cU :=\left(\partial_{x_1}
      \cU(\underline x),\dots, \partial_{x_{e(\Gamma)}}
      \cU(\underline x)\right)$.
    \item[\textbf{Step~3:}] Thanks to steps~1 and~2 the differential form
      \begin{equation}\label{e:betadef}
  \beta^{\underline a}=  \sum_{1\leq i<j\leq e(\Gamma)} {x_i
    C^{\underline a}_j  (\underline x)-x_j
   C^{\underline a}_i(\underline x)\over
  \cF^{\mathsf{a}-1}}\,\Omega_\Gamma^{\epsilon,\kappa}\, 
 dx_1\wedge \cdots \wedge \widehat{dx_i}\wedge \cdots \wedge\widehat{dx_j}\wedge
  \cdots \wedge dx_{e(\Gamma)} \, .
\end{equation}
satisfies the property
  \begin{equation}
     d\beta_\Gamma^{\underline a}= (\mathsf{a}-1) {  P^{\underline
         a}(\underline x)\over
     \cF^{\mathsf{a}}}\Omega_\Gamma^{\epsilon,\kappa}+ {\vec\nabla \cdot \vec C^{\underline a}
     (\underline x)
+ \lambda_U\, 
    c^{\underline a}(\underline x)\over \cF^{\mathsf{a}-1}}\,\Omega_\Gamma^{\epsilon,\kappa} \,.
\end{equation}
where
$\lambda_U=e(\Gamma)-(L+1)(\delta-\epsilon)+\kappa\sum_{i=1}^{e(\Gamma)}
\mu_i$  is the power of $\cU$ in $\Omega_\Gamma^{\epsilon,\kappa}$.
   \end{itemize}

   \medskip
   Therefore, for a given set of derivatives, we have performed the pole
   reduction
   \begin{multline}\label{e:PDEred}
\sum_{\mathsf{a}=a_1+\cdots+a_r\atop
  a_i\geq0} \!\!\!\!\!\!\! c_{\underline a}(\vec m,\vec s;\epsilon,\kappa)\left(\partial\over \partial z_1\right)^{a_1}\cdots\left(\partial
  \over\partial z_r\right)^{a_r} \Omega_\Gamma^{\epsilon,\kappa}\cr
= \sum_{\mathsf{a}=a_1+\cdots+a_r\atop a_i\geq0}
  \!\!\!\!\!\!\!  {\vec\nabla \cdot \vec C^{\underline a}
     (\underline x)
+ \lambda_U\, 
    c^{\underline a}(\underline x)\over  (\mathsf{a}-1) \cF^{\mathsf{a}-1}}\,\Omega_\Gamma^{\epsilon,\kappa}+d\beta^{\epsilon,\kappa}_\Gamma,
\end{multline}
Iterating this procedure gives the partial differential equation~\eqref{e:PFOmegaGeneric}.
 We refer to~\cite{delaCruz:2024xit} for a proof of this reduction and
 details.

\begin{remark}
We remark that this way of solving the linear system includes implicitly the
freedom given by the syzygies of $\textrm{Jac}(\cF)$ and
$\textrm{Jac}(\cU)$ since they belong to
the kernel of the linear systems from~\eqref{e:RedF} and~\eqref{e:RedU}
respectively. It was noticed in~\cite{Lairez:2022zkj},  that in
  the rational case, only the first order syzygies are needed to take
  into account the non-isolated singularities of Feynman integrals.
\end{remark}

\part{Differential operators for various graphs}\label{sec:diffop}

In this section we present some  differential equations acting
on the regulated Feynman integrals. The twist does not change
the singular locus of the integrand. Consequently, the real singularities of the
associated differential operator are the same as when there are no twist
($\epsilon=\kappa=0$). We illustrate the effects of the twist one some
classes of period integrals.

\section{Hypergeometric differential operators: the
                massless box graph}\label{sec:hyper}

\begin{figure}[ht]
	\centering
	\begin{tikzpicture}[scale=0.7]
		\draw[dashed](-1,-1)--(-1.3,-1.3);
		\draw[dashed](1,-1)--(1.3,-1.3);
		\draw[dashed](1,1)--(1.3,1.3);
		\draw[dashed](-1,1)--(-1.3,1.3);
		\draw[dashed, thick](-1,-1)--(-1.0,1.0)--(1.0,1.0)--(1.0,-1.0)--cycle;
		\node[text width=0.5cm, text centered ] at (-1.5,1.5){$p_1$};
		\node[text width=0.5cm, text centered ] at (1.5,1.5){$p_2$};
		\node[text width=0.5cm, text centered ] at (1.5,-1.5){$p_3$};
		\node[text width=0.5cm, text centered ] at (-1.5,-1.5){$p_4$};
	\end{tikzpicture} 
\caption{The box  graph with massless external and internal states. }\label{fig:box}
      \end{figure}
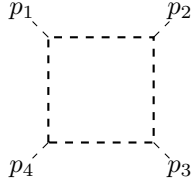

For the box graph in Fig.~\ref{fig:box}, we have four massless momenta
$p_i$ with $1\leq i\leq4$ such that $p_1+\cdots+p_4=0$ and
$p_1^2=\cdots =p_4^2=0$. After scaling the integral and setting $X=
{p_1\cdot p_4\over p_1\cdot p_2}$ the graph
      polynomials are given by
\begin{equation}
	\label{e:Boxgraphpolynomials}
	\cU_{\Box}=x_1+\cdots +x_4,\qquad
	\cF_{\Box}(X)= x_2x_4+X x_1x_3 \, .
      \end{equation}
      The dimensionally regulated       Feynman integral in $D=4-2\epsilon$ is
      given by the       twisted differential 
\begin{equation}\label{e:OmegaBox}
	I_\Box^{\epsilon,0}(X)=\int_{\Delta_4} \left(\cU_{\Box}^2\over \cF_{\Box}(X)\right)^\epsilon
      { \Omega_0\over \cF_{\Box}(X)^2
	}\, ,
      \end{equation}
      with $\Omega_0=x_1 dx_2\wedge dx_3\wedge dx_4-x_2 dx_1\wedge
      dx_3\wedge dx_4+x_3 dx_1\wedge dx_2\wedge dx_4-x_4 dx_1\wedge
      dx_2\wedge dx_3$, 
      and $\Delta_4=\{x_i\geq0, \ 1\leq i\leq 4\}$.
The application of the algorithm described in
section~\ref{sec:griff-dwork-reduct} gives the differential equation
\begin{equation}
    \mathscr{L}_\Box^\epsilon I_\Box^{\epsilon,0}(X)=\frac{(\epsilon +1)\Gamma (-\epsilon -1)^2 }{\Gamma
		(-2 \epsilon )}  \, \left(1+X^{-\epsilon -1}\right) \,,
\end{equation}
with the differential operator
\begin{equation}\label{e:PFBox}
	\mathscr{L}_\Box^\epsilon= (X+1) X {d\over dX}+1+X+\epsilon  \,.
      \end{equation}
      The Feynman integral integrates to hypergeometric functions,
      with the $\epsilon$ expansion as polylogarithms because we have
a variation of      mixed Tate motives~\cite{Goncharov:2001iea,Brown:2013gia}
      \begin{multline}
                	I_\Box^{\epsilon,0}(X)=\frac{4}{\epsilon^{2}
                          X}-\frac{4+2 \ln \! \left(X \right)}{X \epsilon}+\frac{2 \ln \! \left(X \right)-\frac{5 \pi^{2}}{3}+4}{X}\cr
                        +\frac{\epsilon}{X}\,\Big(2 \,\mathrm{Li}_{3}\! \left(-X \right) -2 \ln \! \left(X \right)
                        \mathrm{Li}_{2}\! \left(-X \right)  -\ln \! \left(X \right)^{2} \ln \! \left(X +1\right)+\frac{\ln \left(X
                          \right)^{3}}{3}+\frac{4 \ln \left(X \right)
                          \pi^{2}}{3}\cr-2 \ln \! \left(X \right)-10
                        \zeta \! \left(3\right)+\frac{5
                          \pi^{2}}{3}-4-\pi^{2} \ln \! \left(X
                          +1\right)
                       \Big)+O(\epsilon^2)\,,
      \end{multline}
where $\mathrm{Li}_r(X)=\sum_{n\geq1} X^n/n^r$ are the polylogarithms.      
\section{Hyperelliptic differential operators: Planar two-loop graphs}
\label{sec:motives-two-loop}

\begin{figure}[ht]
\begin{tikzpicture}[scale=0.6]
\filldraw [color = black, fill=none, very thick] (0,0) circle (2cm);
\draw [black,very thick] (-2,0) to (2,0);
\draw [black,very thick] (-2,0) to (-3,0);
\draw [black,very thick] (2,0) to (3,0);
\draw [black,very thick] (0,2) to (0,3);
\draw [black,very thick] (1.414,1.414) to (2.25,2.25);
\draw [black,very thick] (-1.414,1.414) to (-2.25,2.25);
\draw [black,very thick] (0,-2) to (0,-3);
\end{tikzpicture}
\caption{A two-loop graphs with $a=4$, $b=1$ and $c=2$.}\label{fig:a1cgraphs}
\end{figure}
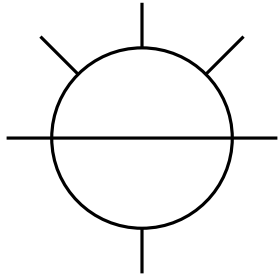

Two-loop graphs can be labelled by the number of edges $(a,b,c)$ on each cycle,
in figure~\ref{fig:a1cgraphs} we have represented a graph with $a=4$, $b=1$, $c=2$.
Planar graphs are graph for which least one edge number is equal to 1,
and their  graph polynomials are given by

\begin{align}\label{eq:UVF}
{\cU}_{(a,1,c)} & = \left(z +  \sum_{i=1}^c x_i\right) \left( \sum_{i=1}^a y_i \right) + z \left( \sum_{i=1}^c x_i\right) ,\\
{\cV}_{(a,1,c);D} & = z\left( \sum_{i=1}^c\sum_{j=1}^a r_{ij}^2 x_i
                      y_j \right)+\left( z + \sum_{i=1}^a
                  y_i\right)\left(\sum_{1\leq i<j\leq c} p_{ij}^2 x_i
                    x_j\right) \cr
                    &+ \left(z+\sum_{i=1}^c x_i\right)\left(
                  \sum_{1\leq i<j\leq a}
                  q_{ij}^2 y_i y_j \right),\cr
\nonumber {\cF}_{(a,1,c);D} & = {\cU}_{(a,1,c)}\left( \sum_{i=1}^c m_{i+a}^2 x_i + \sum_{i=1}^a m^2_i y_i + m^2_{a+c+1} z\right) - {\cV}_{(a,1,c);D}.
\end{align}
\medskip
We define the class of mixed Hodge structures (MHS) than can arise from the
two-loop Feynman integrals:
\begin{definition}[MHS for planar two-loop graphs]\hfill \break\vspace{-.5cm}
\begin{enumerate}[(1)]
\item Let ${\bf MHS}_\mathbb{Q}$ denote the Abelian category of $\mathbb{Q}$-mixed Hodge structures. 

\item The largest extension-closed subcategory of ${\bf MHS}_\mathbb{Q}$ containing the Tate twists of $\mathrm{H}^1(C;\mathbb{Q})$ for every {hyperelliptic curve} $C$ is called ${\bf MHS}_\mathbb{Q}^\mathrm{hyp}$.

\item The largest extension-closed subcategory of ${\bf MHS}_\mathbb{Q}$ containing the Tate twists of $\mathrm{H}^1(E;\mathbb{Q})$ for every {elliptic curve} $E$ is called ${\bf MHS}_\mathbb{Q}^\mathrm{ell}$.
\end{enumerate}
\end{definition}

The main theorem of~\cite{Doran:2023yzu} states for generic values of
$(a,1,c)$ and the spacetime dimension $D$ the Feynman integral attached the
planar two-loop graphs
\begin{equation}
  I_{(a,1,c)}^{D,\underline \nu}= \int_{\Delta_{a+1+c}}
  {\cU^{\sum_{i=1}^{a+1+c} \nu_i-{3D\over2}}\over
    \cF^{\sum_{i=1}^{a+1+c} \nu_i-D}} \prod_{i=1}^a x_i^{\nu_i}
  \prod_{i=1}^c y_i^{\nu_{a+i}}  z^{\nu_{a+c+1}} \Omega_0\,.
\end{equation}
 is a (relative) period integrals of   ${\bf
   MHS}_\mathbb{Q}^\mathrm{hyp}$ because the singular locus is determined
 by the vanishing locus of $\cF_{(a,1,c);D}$
\begin{theorem}[DHV~\cite{Doran:2023yzu}]\label{thm:sunsetmot}
For any values of $a,c$, the cohomology groups of
$X_{(a,1,c);D}=\{\cF_{(a,1,c);D}=0| (x_i,y_i,z)\in\mathbb P^{a+c}\}$ are contained in ${\bf MHS}_\mathbb{Q}^\mathrm{hyp}$.
\end{theorem}
Depending on the value of parameters $D$, $a$, $c$ and the kinematic
invariants $p_{ij}^2$, $q_{ij}^2$, $r_{ij}^2$, and the internal masses
$m_i^2$, the singularities of the integrand change and Feynman
integral can be become a period integral of ${\bf
  MHS}_\mathbb{Q}^\mathrm{ell}$ or just ${\bf MHS}_\mathbb{Q}$.

\subsection{The elliptic curve: the two-loop sunset graph}
\label{sec:elliptic-curve}

We present the case of the sunset integral attached to the two-loop
graph with $a=b=c=1$.

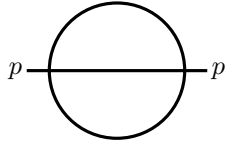
\begin{figure}[!h]
	\centering
	\begin{tikzpicture}[scale=0.3]
		\draw[very thick] (0,0) ellipse (3cm and 3cm);
		\draw[very thick] (-4,0)--(4,0);
			\node [text width=0.5cm, text centered ] at (-4.5,0)
                {$p$};
                	\node [text width=0.5cm, text centered ] at (4.5,0) {$p$};
	\end{tikzpicture} 
	\caption{Two-loop sunset with $a=b=c=1$}
  \label{fig:sunset}
\end{figure}
The graph polynomials are, in the notations introduced above,
\begin{align}\label{eq:UVF111}
{\cU}_{(1,1,1)} & = x_1z+x_1y_1+ y_1z ,\\
{\cV}_{(1,1,1);D} & =p^2 x_1y_1z,\cr
\nonumber {\cF}_{(1,1,1);D} & = {\cU}_{(1,1,1)}\times\left( m_{1}^2 x_1 +
                              m^2_2 y_2 + m^2_{3} z\right) - {\cV}_{(1,1,1);D}.
\end{align}
and the Feynman integral for $D=2-2\epsilon$  and
$\nu_1=\nu_2=\nu_3=1$ reads
\begin{equation}\label{e:OmegaSunsetTwisted}
  I_{(1,1,1)}^{\epsilon}= \int_{\Delta_3}
  \left(\cU_{(1,1,1)}^3\over \cF_{(1,1,1)}^2\right)^\epsilon  \,
  {z dx_1\wedge dy_1-y_1 dx_1\wedge dz+x_1 dy_1\wedge dz\over \cF_{(1,1,1)}}.
\end{equation}
For $\epsilon=0$ it  is shown  in~\cite{Bloch:2013tra,Bloch:2016izu}
that the integral is a regulator period integral of the  ${\bf
  MHS}_\mathbb{Q}^\mathrm{ell}$  associated to the  elliptic curve
defined by $\cF_{(1,1,1)}=0$ in $\mathbb P^2$.

\medskip
\noindent\textbf{The equal-mass case:} When all the mass parameters
are the same $m_1=m_2=m_3$ and setting $t=p^2/m_1^2$, the sunset Feynman  
satisfies the differential equation
\begin{equation}
  \left(   \mathscr{L}^{(2)}_{(1,1,1)} +\epsilon
    \mathscr{L}^{(1)}_{(1,1,1)}  +\epsilon^2  \mathscr{L}^{(0)}_{(1,1,1)} \right)I_{(1,1,1)}^{\epsilon}=-6
     {\Gamma(1+\epsilon)^2\over \Gamma(1+2\epsilon)}
\end{equation}
with the differential operators $\mathscr{L}^{(r)}_{\su(3)} $ of order $r$
\begin{align}
  \label{e:PF2sunset1massepsilon}
     \mathscr{L}^{(2)}_{(1,1,1)} &={d\over dt}\left( t(t-1)(t-9)
                                {d\over dt}\right)+(t-3),\cr
                        \mathscr{L}^{(1)}_{(1,1,1)} &=          (3t^2-10t-9) {d\over
         dt}+3t-5,\cr
    \mathscr{L}^{(0)}_{(1,1,1)} &=       2(t+1).
   \end{align}
  The $\epsilon=0$ piece is the Picard-Fuchs operator associated to the
  modular curve $X_1(6)$ defined by
  $(x_1y_1+x_1y+y_1z)(x_1+y_1+z)=tx_1y_1z$~\cite{Bloch:2013tra}. The
  twist induces the differential $\mathscr{L}^{1}_{(1,1,1)}$ and
  $\mathscr{L}^{0}_{(1,1,1)}$ without affecting the real singularities of
  the total differential operator.

  \medskip

\noindent\textbf{The different mass case:}     For the non-equal-mass case $m_1\neq m_2\neq m_3$ the order of the differential equation is
four with the following $\epsilon$
expansion~\cite{Remiddi:2013joa,Remiddi:2016gno,delaCruz:2024xit}
\begin{equation}\label{e:PF3massSunsetEpsilon}
  \underbrace{\left(  \mathscr{L}^{(4)}_{0}+\sum_{r=0}^4\epsilon^{1+r}
    \hat{\mathscr{L}}^{(4-r)}_{r+1}\right)}_{=:   \mathscr{L}^{\epsilon}_{(1,1,1)}}   I_{(1,1,1)}^{\epsilon}= \mathscr{S}_{(1,1,1)}^\epsilon\,.
\end{equation}
where the differential operators $
\mathscr{L}^{(r)}_{s}$ 
are of order $r$.
The $\epsilon=0$ piece is a fourth order differential operator that factorizes
\begin{equation}\label{e:LsunsetFac}
     \mathscr{L}^{(4)}_{0} =   \mathscr{L}^{(1)}_a\circ
                                \mathscr{L}^{(1)}_b\circ
                                \mathscr{L}^{3-\rm mass}_{\su(3)} 
                              \end{equation}
where     $ \mathscr{L}^{(1)}_a$ and $  \mathscr{L}^{(1)}_b$ are order
one differential operators.
   The differential operator $\mathscr{L}^{3-\rm
     mass}_{\su(3)}$ is the Picard-Fuchs operator 
  for the three masses two-loop sunset integral in two
   dimensions determined by the sunset elliptic curve
   $(x_1y_1+x_1z+y_1z)(m_1^2x_1+m_2^2y_1+m_3^2z)=p^2x_1y_1z$~\cite{Bloch:2016izu}.

The differential operators $ \hat{\mathscr{L}}^{(4-r)}_{r+1}$  are irreducible differential operators
   of  order $4-r$,  therefore the differential
   operator~\eqref{e:PF3massSunsetEpsilon} is irreducible for generic
   values of $\epsilon$.
   
\medskip
   
   The factorisation of the differential operator~\eqref{e:LsunsetFac}
   is understood from the fact that the
   Feynman integral are (relative) periods of the motive  ${\bf
  MHS}_\mathbb{Q}^\mathrm{ell}$ attached to the sunset elliptic curve.
As shown in~\cite{Doran:2023yzu} such factorisation appears for the
$(a,1,c)$ graphs:
\begin{theorem}[Factorisation of differential operators]\label{p:fact-diff}
For any $a,c$, the operator $\mathscr{L}_{(a,1,c);D}$ admits a factorisation
\[
\mathscr{L}_1\mathscr{L}_2 \dots \mathscr{L}_k
\]
where $\Sol(\mathscr{L}_i)$ is either:
\begin{enumerate}[(a)]
    \item a local system with finite order monodromy or
    \item a subquotient of the local system underlying a family of hyperelliptic curves over a Zariski open subset of $\mathbb{A}^1$.
\end{enumerate}
In particular, if $a$ or $c$ is $\leq 2$ then the monodromy representation of $\Sol(\mathscr{L}_i)$ is either 
\begin{enumerate}[(a)]
    \item finite, or
    \item a finite index subgroup of $\mathrm{SL}_2(\mathbb{Z})$. 
\end{enumerate}
\end{theorem}

   \medskip

The twist in the differential form~\eqref{e:OmegaSunsetTwisted} induces an $\epsilon$ deformation
of the Picard-Fuchs operator $\mathscr{L}^{\epsilon}_{(1,1,1)}$. This
does not affect the real singularities
of the differential operator because the $\epsilon$ factor
in~\eqref{e:OmegaSunsetTwisted} does not change the nature of the singular
locus which is still given by the same elliptic curve as in the
$\epsilon=0$ case. Therefore, the $\epsilon$ deformation only affects
the local monodromies and the apparent
singularities of the differential operator,  as can be seen from the
coefficient of the highest order term
\begin{multline}
  \mathscr{L}^{\epsilon}_{(1,1,1)} \Big\vert_{(d/dt)^4}= (p^2)^3\prod_{i=1}^4 (p^2-
  \mu_i^2) \Big(-\left(2 \epsilon +5\right) (p^2)^{2}-2
    \left(m_{1}^{2}+m_{2}^{2}+m_{3}^{2}\right) \left(1+2 \epsilon
    \right)p^2\cr+ \left(7+6 \epsilon \right)\prod_{i=1}^4 \mu_i
\Big)  \, , 
\end{multline}
where $\mu_i=\{m_1+m_2+m_3,-m_1+m_2+m_3,m_1-m_2+m_3,m_1+m_2-m_3\}$ are
the thresholds.

\medskip
The action  of $\mathscr{L}^{\epsilon}_{(1,1,1)}$ leads to the
inhomogeneous differential equation~\eqref{e:PF3massSunsetEpsilon},
with an inhomogeneous term given by
   \begin{equation}
     \mathscr{S}_\su(\vec m,t,\epsilon)=\frac{c_{23}(t,\epsilon)\Gamma (\epsilon +1)^2}{ (m_{2} m_{3})^{2 \epsilon}\Gamma (1+2\epsilon)}+\frac{c_{13}(t,\epsilon)\Gamma (\epsilon +1)^2}{ (m_{1} m_{3})^{2 \epsilon }\Gamma (1+2
   \epsilon)}+\frac{c_{12}(t,\epsilon)\Gamma (\epsilon +1)^2}{ (m_{1} m_{2})^{2 \epsilon }\Gamma (1+2
   \epsilon )} \, , 
\end{equation}
 obtained from the integration of the exact differential
 $d\beta_{(1,1,1)}^{\epsilon}$ in~\eqref{e:PFOmegaGeneric}.
 The coefficients 
$c_{12}(t,\epsilon)$, $c_{13}(t,\epsilon)$ and
$c_{23}(t,\epsilon)$ are polynomials of degree 4  in $t$ and degree 2 in
$\epsilon$, respectively which expressions are given on the {\tt
  SageMath} worksheet \href{https://nbviewer.org/github/pierrevanhove/TwistedGriffithsDwork/blob/main/Worksheets/Sunset-Twoloop-3mass-Epsilon.ipynb}{Sunset-Twoloop-3mass-Epsilon.ipynb}.
Expanding in powers of $\epsilon$, we have
\begin{equation}
    \mathscr{S}_\su(\vec m,t,\epsilon)=\mathscr{S}_\su^0(\vec m,t)
    +\left(c^{(1)}_0(\vec m)+
    \sum_{i=1}^3  c^{(1)}_i(\vec m)\log(m_i)\right)\,  \epsilon+O(\epsilon^2)
\end{equation}
with  the leading term given by
   \begin{equation}
   \mathscr{S}_\su^0(\vec m,t)=60 t^{4}+56\left( m_{1}^{2}+ m_{2}^{2}+
     m_{3}^{2}\right) t^{3}
   -308 \prod_{i=1}^4 \mu_i.
 \end{equation}
For $\epsilon=0$ the two-loop sunset integral satisfies  the
differential equation~\cite{Adams:2013nia,Bloch:2016izu}
\begin{equation}
  \mathscr{L}^{3-\rm mass}_\su f_\su^{(0)}(t)= s_0(\vec m,t)+ \sum_{i=1}^3
  s_i(\vec m,t)\log(m_i^2)  \,.
\end{equation}
It can be checked that 
\begin{equation}
  \mathscr{S}_\su^0(\vec m,t)=   \mathscr{L}^{(1)}_1
     \mathscr{L}^{(2)}_1    \mathscr{L}^{3-\rm mass}_\su  f_\su^{(0)}(t),
   \end{equation}
   showing that the structure of the inhomogeneous term is compatible
   with the factorisation of the $\epsilon=0$ piece of the
   differential operator in~\eqref{e:PF3massSunsetEpsilon}.

\section{Calabi--Yau differential operators: sunset multiloop graphs}
\label{sec:CY}

\begin{figure}[!h]
	\centering
	\begin{tikzpicture}[scale=0.4]
		\draw[very thick] (0,0) ellipse (3cm and 3cm);
		\draw[very thick] (0,0) ellipse (3cm and 2.cm);
		\draw [very thick](0,0) ellipse (3cm and 1.5cm);
		\draw [very thick](0,0) ellipse (3cm and 1.cm);
		\draw[very thick] (-4,0)--(4,0);
		\node[text width=0.5cm, text centered ] at (0,-.3) {$\vdots$};
		\node[text width=0.5cm, text centered ] at (0,0.6) {$\vdots$};
		\node [text width=0.5cm, text centered ] at (-4.5,0)
                {$p$};
                	\node [text width=0.5cm, text centered ] at (4.5,0) {$p$};
	\end{tikzpicture} 
	\caption{Multi-loop sunset with $n$ edges}
  \label{fig:nsunset}
\end{figure}
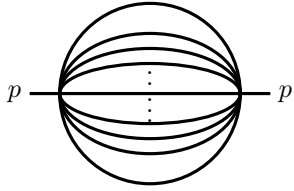

We now turn to the  $n-1$-loop sunset integral in $D=2-2\epsilon$
dimensions attached to the graph in fig.~\ref{fig:nsunset} which reads
\begin{equation}\label{e:IsunsetEpsilon}
  I^{\epsilon}_{\su(n)}(p^2,\vec m,t)= \int_{\Delta_n}\,\left(\cU_{\su(n)}^{n}\over
    \cF_{\su(n)}^{n-1}\right)^\epsilon\, {\sum_{i=1}^n(-1)^{i-1}
    \bigwedge_{j=1\atop j\neq i}^n dx_j\over
    \cF_{\su(n)}},
\end{equation}
with the domain of integration $\Delta_n=\{x_i\geq0, 1\leq i\leq n\}$
and the graph polynomials
\begin{align}\label{e:UFsunset}
  \cU_{\su(n)}&= x_1\cdots x_{n}\sum_{i=1}^n {1\over x_i}\,, \cr
     \cF_{\su(n)}&= \cU_{\su(n)} \times\sum_{i=1}^n m_i^2x_i-p^2 x_1\cdots x_n \, .
\end{align}
Notice that $\cU_{\su(n)}^n/ \cF_{\su(n)}^{n-1}$ is a homogeneous
rational function  of degree
0 in $(x_1,\dots,x_n)$. As usual the differential form is defined in
the complement of the vanishing locus of the denominator in $X_n=\{\cF_{\su(n)}=0\}$.
The Feynman integral
being a (relative) period of a Calabi--Yau  manifold of complex
dimension $n-2$ defined by the equation $\cF_{\su(n)}=0$~\cite{Bloch:2013tra,Bloch:2014qca,Bloch:2016izu,Bourjaily:2019hmc,Bonisch:2020qmm,Bonisch:2021yfw,Candelas:2021lkc,Forum:2022lpz}.

\bigskip
\noindent\textbf{The equal-mass case:}
For the equal-mass case $m_1=\cdots =m_{n}$ and $p^2=t\,m_{n}^2$ in~\eqref{e:UFsunset} the sunset
Feynman integral satisfies the differential equation
\begin{equation}
\underbrace{\left( \sum_{r=0}^{n-1} \epsilon^r
    \mathscr{L}_{\su(n)}^{(n-1-r)}\right)}_{=: \mathscr{L}_{\su(n)}^{\epsilon}} I_\su^\epsilon(t)=
  -n! {\Gamma(1+\epsilon)^{n-1}\over \Gamma(1+(n-1)\epsilon)}.
\end{equation}
Like the case $n=3$ described above the term of order $\epsilon^r$ is
a differential operator $\mathscr{L}_{\su(n)}^{(r)} $ of order $r$ in $t$.
The coefficient of $\epsilon^0$ is the differential operator of order
$n-1$  derived
  in~\cite{Vanhove:2014wqa}  (see as
  well~\cite{Bonisch:2020qmm,Pogel:2022yat,Pogel:2022ken,Pogel:2022vat,Mishnyakov:2023wpd,Mishnyakov:2023sly}).

  The $\epsilon$ deformation does not change the real singularities of
the differential operators because the twist
in~\eqref{e:IsunsetEpsilon} does not introduce new singularities.

  \medskip
This is seen as well on the form of the $\epsilon$-deformed
Picard-Fuchs operator for the $n=4$ case with $m_1=m_2=m_3=m_4$ 
  \begin{multline}
    \mathscr{L}_{\su(4)}^{\epsilon}=
    -(t - 16)  (t - 4)  t^2\left(d\over dt\right)^3 - 6  (t^3 -
                              15  t^2 + 32  t)  \left(d\over dt\right)^2 - (7  t^2 - 68  t +
                              64)  \left(d\over dt\right) - t + 4\cr
                              +\epsilon \left(-6  (t - 10  ) t^2 \left(d\over dt\right)^2 -
      6  (3  t - 20  ) t \left(d\over dt\right) +18- 6  t 
      \right)\cr
    +\epsilon^2\left(-(11  t^2 -
      28  t - 64)  \left(d\over dt\right) - 11  t + 14\right)+ \epsilon^3\left(-6  t - 12\right) \,,
  \end{multline}
 where the $\mathscr{L}_{\su(4)}^{\epsilon}|_{\epsilon=0}$ operator is the
  Picard-Fuchs operator for the $K3$ surface with Picard number 19~\cite{Bloch:2014qca}.

  \medskip
\noindent\textbf{Different mass case:} We represent the result for the four different masses $m_1\neq m_2\neq m_3
\neq m_4$ for the case $n=4$. In that case the singular locus is a
$K3$ surface of Picard number 16~\cite{Lairez:2022zkj}.
The $\epsilon$-deformed 
 differential operator has order 11 and has the $\epsilon$ expansion 
  \begin{equation}
    \mathscr{L}_{\su(4)}^{\epsilon}=     \sum_{r=0}^{16}\epsilon^r
    \mathscr{L}^{(11)}_{r}+\sum_{r=0}^{11} \epsilon^{16+r}  \mathscr{L}^{(11-r)}_{16+r}\, ,
  \end{equation}
  where the differential operators $\mathscr{L}^{(r)}_s$ are of order
  $r$.  The order 11 part of this deformed Picard-Fuchs operator has
  degree 16 in $\epsilon$.

  \medskip
    The order $\epsilon=0$ operator factorises as
  \begin{equation}
         \mathscr{L}^{(11)}_{0}=   \mathscr{L}^{(1)}_{a_1} \circ
         \cdots  \circ \mathscr{L}^{(1)}_{a_5}  \circ
         \mathscr{L}^{4-\rm mass}_{\su(4)}\,,
       \end{equation}
        where  $ \mathscr{L}^{(1)}_{a_1},\dots,  \mathscr{L}^{(1)}_{a_5}$ are  first order operators
       and $\mathscr{L}^{4-\rm mass}_{\su(4)}$ is the sixth  order differential operator for the
       three-loop sunset integral with the all different mass configurations
       given in section~4.3 of~\cite{Lairez:2022zkj}.

       \medskip

       The coefficient of the highest order term $(d/dt)^{11}$    is given by
       \begin{multline}
                   \mathscr{L}^{\epsilon}_{\su(4)}\Big\vert_{(d/dt)^{11}}=
                   t^{11}\left(t-(m_{1}+m_{2}-m_{3}-m_{4})^2\right) \cr\times
   \left(t-(m_{1}-m_{2}+m_{3}-m_{4})^2\right)
   \left(t-(m_{1}+m_{2}+m_{3}-m_{4})^2\right) \cr\times
   \left(t-(m_{1}-m_{2}-m_{3}+m_{4})^2\right)
   \left(t-(m_{1}+m_{2}-m_{3}+m_{4})^2\right)\cr\times
   \left(t-(m_{1}-m_{2}+m_{3}+m_{4})^2\right)
   \left(t-(-m_{1}+m_{2}+m_{3}+m_{4})^2\right) \cr\times
   \left(t-(m_{1}+m_{2}+m_{3}+m_{4})^2\right)
                   \, q^{[1111]}(t,\epsilon).
                 \end{multline}
                 The $\epsilon$ dependence appears only in the
                 apparent singularities determined by the polynomial
                 $q^{[1111]}(t,\epsilon)$ of degree 17 in
                 $t$ and 16 in $\epsilon$. The polynomial is given
                 in the only worksheet \href{https://nbviewer.org/github/pierrevanhove/TwistedGriffithsDwork/blob/main/Worksheets/Sunset-Threeloop-Epsilon.ipynb}{Sunset-Threeloop-Epsilon.ipynb}.

\section{Discussion}
We have presented a generalisation of the Griffiths-Dwork reduction
for deriving differential operators acting on Feynman integrals in
dimensional or analytic regularisation. The algorithm
makes a special use of the fact that the twist from the regularisations
is the power of a degree zero homogeneous rational function build from
the graph polynomials.

\medskip

The algorithm gives the minimal order  (non-factorisable) D-module of differential
operators acting on regulated Feynman
integrals. At each derivative order the
 procedure consists of solving the linear systems from the reductions
 with the respect the Jacobian ideal of the graph polynomials $\cF$
 in~\eqref{e:RedF} and $\cU$  in~\eqref{e:RedU}  in order to
determine the coefficients $c_{\underline a}(\underline z)$ and the 
inhomogeneous term $\beta^{\underline a}_\Gamma$
in~\eqref{e:PFOmegaGeneric}.

\medskip

Because the twisted differential $\Omega_\Gamma^{\epsilon,\kappa}$
has the same singularities as  $\Omega_\Gamma^{0,0}$, the regularisation parameters $\epsilon$ or
  $\kappa$ do not affect the discriminant locus but only the local
  monodromies.  This reflects on the fact that these parameters only
  affect the apparent singularities of the differential operators.

\medskip
  
With this algorithm we can derive a Gr\"obner basis of partial
differential operators in some multiple scale cases. The differential
operators produced by the algorithm of this paper might arise as
specialisation of the system of partial differential operators
obtained by GKZ  approach. The restriction of the GKZ  D-module is a
difficult open problem, which we leave for further investigations.

\section*{Acknowledgements}
I would like to thank Spencer Bloch, Matt Kerr,
Leonardo de la Cruz, Pierre Lairez, Eric Pichon-Pharabod, Andrew
Harder and Charles Doran for very enjoyable collaboration on the
topics reported on this text.

I would like to thank the organisers of the Regulator V conference
(3-13 June 2024 in Pisa, Italy)  for the opportunity to present these results.

The author was supported in part by French National Agency for Research
grant ``Observables'' (ANR-24-CE31-7996). This research was supported by the Munich Institute for Astro-, Particle and BioPhysics (MIAPbP) which is funded by the Deutsche Forschungsgemeinschaft (DFG, German Research Foundation) under Germany's Excellence Strategy - EXC-2094-390783311.

\end{document}